\documentclass{amsart}

\newtheorem{theorem}{Theorem}[section]
\newtheorem{lemma}[theorem]{Lemma}

\newtheorem{proposition}[theorem]{Proposition}

\theoremstyle{definition}
\newtheorem{definition}[theorem]{Definition}
\newtheorem{example}[theorem]{Example}

\theoremstyle{remark}

\numberwithin{equation}{section}

\begin{document}

\title{On right coideal subalgebras}
\author{V.K. Kharchenko}
\address{Universidad Nacioanal Aut\'onoma de M\'exico, Facultad de Estudios Superiores Cuautitl\'an, Primero de Mayo s/n, Campo 1, CIT, Cuautitlan Izcalli, Edstado de M\'exico, 54768, MEXICO}
\email{vlad@servidor.unam.mx}
\thanks{The author was supported by PAPIIT IN 108306-3, UNAM, M\'exico.}



\begin{abstract}
Let $H$ be a character Hopf algebra. Every  
right coideal subalgebra that contains the coradical has a PBW-basis
which can be extended up to a PBW-basis of $H.$ 
\end{abstract}

\maketitle

\section{Introduction}

In the present paper we consider right coideal subalgebras in character 
Hopf algebras, the Hopf algebras generated by skew-primitive semi-invariants.
This class includes quantum enveloping algebras of Kac--Moody algebras, 
all their generalizations 
(see,  G. Benkart, S.-J. Kang, and D. Melville \cite{BKM}; 
M. Costantini, and M. Varagnolo \cite{CV}; S.-J. Kang \cite{Kan}),
the bosonizations of quantum symmetric algebras
related to diagonal braidings, and so on.  We prove the following  general statement on 
the structure of the right coideal subalgebras.
\begin{theorem} Let $H$ be a character Hopf algebra. Every  
right coideal subalgebra that contains all group-like elements has a PBW-basis
which can be extended up to a PBW-basis of $H.$
\label{mainint}
\end{theorem}
This theorem is new for Hopf subalgebras too.

\section{Preliminaries}

\noindent 
\subsection{PBW-generators}
Let $S$ be an algebra over a field {\bf k} and $A$ is a subspace of $S$
with a basis $\{ a_j\, |\,  j\in J\} .$ A linearly ordered subset $V\subseteq S$ is said to be
a {\it set of PBW-generators of $S$ over $A$} if there exists 
a function $h:V\rightarrow {\bf Z^+}\cup \{ \infty \},$
called the {\it height function}, such that the set of all products
\begin{equation}
a_jv_1^{n_1}v_2^{n_2}\cdots v_k^{n_k}, 
\label{pbge}
\end{equation}
where $j\in J,\ \ v_1<v_2<\ldots <v_k\in V,\ \ n_i<h(v_i), 1\leq i\leq k$
is a basis of $S.$ The value $h(v)$ is referred to as the {\it height} of $v$ in $V.$
If $A={\bf k}$ is the ground field, then we shall call $V$  simply as a {\it set of PBW-generators of} $S.$ 
\begin{definition} \rm
Let $V$ be a set of PBW-generators of $S$ over a subalgebra $A.$
Suppose that the set of all words in $V$ as a free monoid has its own order $\prec $
(that is, $a \prec b$ implies $cad\prec cbd $ for all words $a, b, c, d$ in $V$). 

1) A {\it leading word} of $s\in S$ is the maximal word $m=v_1^{n_1}v_2^{n_2}\cdots v_k^{n_k}$
that appears in the decomposition of $s$ in the basis (\ref{pbge}).

2) A {\it leading term} of $s$ is the sum $am$ of all terms $\alpha _ia_im$
that appear in the decomposition of $s$ in the basis (\ref{pbge}), where
$m$ is the leading word of $s.$

3) The order $\prec $ is {\it compatible with the PBW-decomposition related to} $V$
if the leading word of each product $W=w_1w_2\cdots w_m,$ $w_i\in V\cup \{ a_j\}$
(considered as an element of $S$)
 is less than or equal to the word $W^{t}$ that appears from  
$w_1w_2\cdots w_m$ by deletion of all letters $a_j.$ 
\label{pbgd}
\end{definition}

\noindent 
\subsection{Associated graded algebra} 
Let $\Gamma $ be a completely ordered additive (commutative)
monoid. A $\Gamma $-{\it filtration} on an algebra $S$
is a map $D$ from $S$ to $\Gamma $ extended by a symbol $-\infty $ with the properties
$$
\begin{matrix}
1. & D(x)=-\infty \Longleftrightarrow  x= 0, \hfill & D(1\cdot {\bf k}^*)=0,\hfill \cr
2. & D(x-y)\leq {\rm max}\{ D(x), D(y)\}, \hfill & D(xy)\leq D(x)+D(y).\hfill
\end{matrix}
$$
The element $D(x)$ is the {\it degree } of $x.$ We denote by $S_{\gamma }$
the subspace of all elements of degree $\leq \gamma ,$ while by
$S_{\gamma }^{-}$ the subspace of all elements of degree $<\gamma .$
\begin{example} \rm (Filtration by constitution). Suppose that an algebra $S$
is generated over a subalgebra $A$ by a set $X.$ That is, there exists a homomorphism
$\xi : A\langle X\rangle \rightarrow S,$ where  $A\langle X\rangle$ is the free product 
$A\ast {\bf k}\langle X\rangle .$ A {\it constitution} of a word $u$ in $A \cup X$
is a family of non-negative integers $\{ m_x, x\in X\} $
such that $u$ has $m_x$ occurrences of $x.$
Certainly almost all $m_x$ in the constitution are zero.
We fix an arbitrary complete order, $<,$ on the set $X.$

Let $\Gamma $ be the free additive (commutative) monoid generated by $X.$
The monoid $\Gamma $ is a completely ordered monoid with respect to 
the following order:
\begin{equation}
m_1x_{i_1}+m_2x_{i_2}+\ldots +m_kx_{i_k}>
m_1^{\prime }x_{i_1}+m_2^{\prime }x_{i_2}+\ldots +m_k^{\prime }x_{i_k}
\label{ord}
\end{equation}
if the first from the left nonzero number in
$(m_1-m_1^{\prime}, m_2-m_2^{\prime}, \ldots , m_k-m_k^{\prime})$
is positive, where $x_{i_1}>x_{i_2}>\ldots >x_{i_k}$ in $X.$
We associate a formal degree $D(u)=\sum _{x\in X}m_xx\in \Gamma $
to a word $u$ in $A\cup X,$ where $\{ m_x\, | \, x\in X\}$ is the constitution of $u.$
Respectively, if $f=\sum \alpha _i u_i\in A\langle X\rangle ,$ $0\neq \alpha _i\in {\bf k}$
then 
\begin{equation}
D(f)={\rm max}_i\{ D(u_i)\} .
\label{degr}
\end{equation}
 On $S$ we define a $\Gamma $-filtration
related to $\xi : A\ast {\bf k}\langle X\rangle \rightarrow S$ as follows:
\begin{equation}
D(s)={\rm min }\{ D(f)\, |\, \xi (f)=s\}. 
\label{degr1}
\end{equation}
\label{exa}
\end{example} 

With every $\Gamma $-filtered algebra $S$ a $\Gamma $-graded algebra is associated
in the obvious way. For each $\gamma \in \Gamma ,$ write gr$_{\gamma }\, S$ for the linear 
space $S_{\gamma }/S_{\gamma }^- .$ Then
$$
{\rm gr}\, S=\bigoplus _{\gamma \in \Gamma }{\rm gr}_{\gamma }\, S,
$$
with the product defined by
$$
(x+S^-_{\gamma })(y+S^-_{\delta })=xy+S_{\gamma +\delta }^-.
$$
If the filtration is as in Example \ref{exa}, then the algebra gr$\, S$ is generated
by gr$\, X\cong X$ and gr$\, A=A.$
\begin{lemma} If the associated graded algebra {\rm gr}$\, S$ 
for a $\Gamma $-filtered algebra $S$ has a set of
homogeneous PBW-generators, $V=\{ v_i\}, $ over a subspace $A\subseteq S_{0},$ then
each set of representatives $\hat{v}_i,$ $v_i=\hat{v}_i+S_{D(v_i)}^- ,$ is a set
of PBW-generators of $S$ over $A.$  
\label{gruf}
\end{lemma}
This lemma has been proved by S. Ufer \cite[Proposition 46]{Ufe}
for the case  $\Gamma ={\bf Z}^+.$ In the general case the proof is quite similar.

\begin{lemma} Let {\bf U} be a subalgebra of a  
$\Gamma $-filtered algebra $S.$ If {\rm gr}$\, S$ is a free left $($right$)$
$\Gamma $-graded module over {\rm gr}{\bf U} then $S$ is a free left $($right$)$ module over
 {\bf U}.
\label{gruf1}
\end{lemma}
This is part of the folklore.

\noindent 
\subsection{Thin elements and replacement of basis} Suppose that an algebra $S$
has a set $V$ of PBW-generators over a subalgebra $A,$ and $1\in \{ a_j\} .$
Let $\prec $ be any complete order of the free monoid of words in $V$ which is 
compatible with the PBW-decomposition related to $V.$
An element $c \in S$ is said to be a {\it thin element} if its decomposition
in the basis (\ref{pbge}) has the from
\begin{equation}
c=v^m+\sum \alpha _{ij}a_jW_i, 
\label{thin}
\end{equation}
where $W_i\prec v^m,$ and either $m$ divides the height of $v$ in $V,$ or $h(v)=\infty .$ 

Let $T\subseteq S$ be a some set of thin elements. Suppose that  
for each $v\in V$ there exists at most one element 
$c_v\in T$ with the leading term
of the form $v^m,$ $m\geq 1.$   
One may construct a new set of PBW-generators, $P_T,$ related to $T$ 
in the following way.

If in $T$ there does not exist an element with the leading 
term of the form $v^m,$ $m\geq 1$ we include $v$ in $P_T$ and define the height
$h_T(v)$ related to  $T$ to be equal to $h(v).$

If there exists an element $c_v\in T$ with the leading term $v^m$ and
 $m>1,$ we include in $P_T$ both elements: $v$ and $c_v.$
In this case we  define the height of $v$ in $P_T$  
to be equal to $m,$ while the height
of $c_v$ related to $P_T$  is the quotient $h(v)/m.$

If there exists an element $c_v\in T$ with the leading term $v$ 
(that is, if $m=1$) then we include in $P_T$ just $c_v$ (that accidentally may be equal to $v$).
In this case the height of $c_v$ with respect to $T$  by definition equals $h(v).$

We extend on $P_T$ the  order $<$ in the natural way: 
$c_v<w$ if and only if ($v\leq w$ $\& $ $c_v\neq w$); and $c_v<c_w$ if and only if $v<w.$
In particular $c_v<v,$ provided that $m>1,$ where $m$ is defined by $c_v$ in (\ref{thin}).
\begin{lemma} 
The set $P_T$ is a set of PBW-generators of $S$ over $A.$ 
\label{pbw}
\end{lemma}
\begin{proof} We have to prove that the monotonous restricted 
words (\ref{pbge}) with $v_i\in P_T$ are linearly independent in $S,$ and they span $S.$

By definition if $\theta _1<\theta _2\in P_T,$ $\theta _1=v_1^{m_1}+\cdots ,$
$\theta _2=v_2^{m_2}+\cdots ,$
then either $v_1<v_2,$ or $v_1=v_2=v$ with $\theta _1=c_v=v^m+\cdots ,$
$\theta _2=v.$ Therefore each monotonous word in $P_T$ has a form
\begin{equation}
\theta _1^{n_1}v_1^{r_1}\theta _2^{n_2}v_2^{r_2}\cdots \theta _k^{n_k}v_k^{r_k}, 
\label{slo}
\end{equation}
where $v_1<v_2<\cdots <v_k,$ $\theta _i=c_{v_i}=v_i^{m_i}+\cdots .$
Of course if $m_i=1,$ then $r_i=0,$ while if $m_i>1,$ then $r_i\geq 0.$
The word (\ref{slo}) is restricted if and only if 
$n_i<h_T(\theta _i)=h(v_i)$ in the case $m_i=1,$ and 
$n_i<h_T(\theta _i)={h(v_i)/ m_i},$ $r_i<m_i$ otherwise.
If we replace $\theta _i:=v_i^{m_i}+\cdots $
in (\ref{slo}) and then develop multiplication,
we get a linear combination $W+\sum_q \alpha _q W_q$ of words in $V\cup \{ a_j \},$ where 
\begin{equation}
W=v_1^{n_1m_1+r_1}v_2^{n_1m_2+r_2}\cdots v_k^{n_km_k+r_k}.
\label{slo1}
\end{equation} 
Let us as above denote by $W^t$ a word in $V$ that appears from $W$
by deletion of all $a_j.$ 
Since $\prec $ is a monoidal order, we have $W_q^{t}\prec W,$ for all $q.$
If $m_i>1$ then 
$$
n_im_i+r_i\leq (\frac{h(v_i)}{m_i}-1)m_i+(m_i-1)=h(v_i)-1.
$$
If $m_i=1,$ again $n_im_i+r_i=n_i<h(v_i).$
Hence (\ref{slo1}) is a monotonous restricted word, provided that so is (\ref{slo}).
Since the order $\prec $ is compatible with the PBW-decomposition
related to $V,$ all words  $v_1v_2\cdots v_k$ that appears in the PBW-decomposition 
(\ref{pbge}) of $\sum_q \alpha _qW_q$ are less than (\ref{slo1}). 
It is important to note that in this way to different monotonous restricted words in $P_T$ 
correspond different monotonous restricted words in $V$. 

Suppose that a linear combination, $\Xi =\sum \alpha _{ij}  a_jU_i,$ over $A$
of monotonous restricted words
in $P_T$  equals zero in $S.$ 
Let us first in $\Xi $ replace each $\theta \in P_T$ by its representation,
$\theta :=v^m+\cdots ,$ in terms of PBW-generators $V,$ 
and next develop the multiplication.
We get a linear combination $\Xi ^{\prime }=\sum \beta  _kW_k$ of words in 
$P_T\cup \{ a_j \, |\, j\in J  \} .$ Since $\prec $ is a monoidal order
($a\prec b$ implies $cad\prec cbd$), the maximal word $W_0$ among $W_k^t$
appears in the decomposition of summands $a_j U$ with the only word $U.$
If $U$ has form (\ref{slo}), then $W_0$ has form (\ref{slo1}). 
Since the order $\prec $ is compatible with the PBW-decomposition 
related to $V,$ the maximal word of $\Xi ^{\prime }$ 
in  the PBW-decomposition appears just in summands $\alpha_{sj} a_jW_0,$
where the equality $U_s=U$ defines the index $s.$ Hence $\Xi ^{\prime }\neq 0$ in $S.$ A contradiction.

To see that monotonous restricted words in $P_T\cup \{ a_j \}$ 
span $S$ we may use a standard induction on words ordered by $\prec $. Indeed, for any
monotonous restricted word in $V,$
\begin{equation}
W=v_1^{s_1}v_2^{s_2}\cdots v_k^{s_k},
\label{slo2}
\end{equation}
we have $s_i=n_im_i+r_i,$ $r_i<m_i,$ $n_i<h_T(\theta _i),$ where either
$\theta _i=c_{v_i}=v_i^{m_i}+\cdots \in T,$
or the set $T$ has no elements with the leading word $v_i^{m},$ $m\geq 1$
(and hence $\theta _i=v_i\in P_T, m_i=1$).
Suppose by induction that values of all super-words
smaller then $W$ belong to the linear space span by monotonous
restricted words in $P_T\cup \{ a_j \}.$ 
The difference between $W$ and (\ref{slo}) is a linear 
combination of words $U$ in $P_T\cup \{ a_j \} $ such that $U^{t}$ are less than $W.$
Due to the comparability of $\prec $ with the PBW-decomposition
related to $V,$ the PBW-decomposition of any such a word
has only summands $\alpha a_jW^{\prime }$ with $W^{\prime }\prec W.$ 
By induction the lemma is proved. \end{proof}
  
{\bf Remark}. In the proof we have seen that the leading term of (\ref{slo}) in basis (\ref{pbge})
equals (\ref{slo1}). 

\section{PBW-basis of a character Hopf algebra}

Recall that a Hopf algebra $H$ is referred to as a {\it character} Hopf 
algebra if the group $G$ of all grouplike elements is commutative
and $H$ is generated over {\bf k}$[G]$ by skew primitive 
semi-invariants $s_i,\ i\in I:$
\begin{equation}
\Delta (s_i)=s_i\otimes 1+g_{s_i}\otimes s_i,\ \ \ \
g^{-1}s_ig=\chi ^{s_i}(g)s_i, \ \ g, g_{s_i}\in G,
\label{AIc}
\end{equation}
where $\chi ^{s_i},\, i\in I$ are characters of the group $G.$
Let us associate a variable $x_i$ to $s_i.$
For each word $u$
in $X=\{ x_i\, |\, i\in I\}$
we denote by $g_u$
an element of $G$
that appears from $u$
by replacing each $x_i$ with $g_{s_i}.$
In the same way we denote by $\chi ^u$
a character that appears from $u$
by replacing of each $x_i$ with $\chi ^{s_i}.$
We define a bilinear skew commutator by the formula
\begin{equation}
[u,v]=uv-p_{uv}vu,
\label{sqo}
\end{equation}
where by definition  $p_{uv}=\chi ^u(g_v)=p(u,v).$

The group  $G$ acts on the free algebra ${\bf k}\langle X\rangle $
by $ g^{-1}ug=\chi ^u(g)u,$ where $u$ is an arbitrary monomial 
in $X.$
The skew group algebra $G\langle X\rangle $
has the natural Hopf algebra structure 
$$
\Delta (x_i)=x_i\otimes 1+g_{s_i}\otimes x_i, 
\ \ \ i\in I, \ \ \Delta (g)=g\otimes g, \ g\in G.
$$

The construction of the PBW-basis given in
\cite{Kh3} requires $I$ to be finite.
In order to get a PBW-basis in the general case
one may slightly modify this construction as follows.

We fix a Hopf algebra homomorphism
$\xi :G\langle X\rangle \rightarrow H,$ $\xi (x_i)=s_i,$
$\xi (g)=g,$ $i\in I,$ $g\in G,$ and consider the filtration
$D$ related to $\xi $ as defined in Example \ref{exa}.
This filtration is compatible with the Hopf algebra structure, that is 
\begin{equation}
\Delta (H_{\gamma })\subseteq \sum_{\delta +\varepsilon =\gamma} H_{\delta }\otimes H_{\varepsilon},
\ \sigma (H_{\gamma })\subseteq H_{\gamma },
\label{comj}
\end{equation}
where $\sigma $ is the antipode. Therefore gr$\, H$ is also a character Hopf algebra
generated by {\bf k}$[G]$ and $\xi (X).$ 

\underline{
By Lemmas \ref{gruf} and \ref{gruf1} in what follows 
we may suppose 
that $H$ is a $\Gamma $-graded}  \newline \underline{character Hopf algebra} (or, in the other
words, it is homogeneous in each of the generators $s_i$).

On the set of all words in $X$ we fix the lexicographical order
with the priority from the left to the right,
where a proper beginning of a word is considered to 
be greater than the word itself.
This order is not monoidal: we have $x>x^2,$ while $xy<x^2y$ if $y<x.$
It has neither ACC, nor DCC: $x>x^2>\ldots >x^n>\ldots ;$
$xy<x^2y<\ldots x^ny<\ldots .$ By these reasons
we need one more order, the {\it Hall} (or {\it deg-lex}) one,
$u\prec v$ if $D(u)<D(v)$ or $D(u)=D(v)$ $\& $ $u<v.$
Since there exist just a finite number of words in $X$ of a given constitution,
the Hall order indeed is a complete order of the free monoid.  

A non-empty word $u$
is called a {\it standard} word (or {\it Lyndon} word, or 
{\it Lyndon-Shirshov} word) if $vw>wv$
for each decomposition $u=vw$ with non-empty $v,w.$
A {\it nonassociative} word is a word where brackets 
$[, ]$ somehow arranged to show how multiplication applies.
If $[u]$ denotes a nonassociative word then by $u$ we denote 
an associative word obtained from $[u]$ by removing the brackets
(of course $[u]$ is not uniquely defined by $u$ in general).

The set of {\it standard nonassociative} words is the biggest set $SL$
that contains all variables $x_i$
and satisfies the following properties.

1)\ If $[u]=[[v][w]]\in SL$
then $[v],[w]\in SL,$
and $v>w$
are standard.

2)\  If $[u]=[\, [[v_1][v_2]]\, [w]\, ]\in SL$ then $v_2\leq w.$

\noindent
Every standard word has
the only alignment of brackets such that the appeared
nonassociative word is standard (the Shirshov theorem).
In order to find this alignment one may use the following
procedure: The factors $v, w$
of the nonassociative decomposition $[u]=[[v][w]]$
are the standard words such that $u=vw$
and  $v$ has the minimal length (\cite{pSh2}, see also \cite{Lot}).
\begin{definition} \rm A {\it super-letter}
is a polynomial that equals a nonassociative standard word
where the brackets mean (\ref{sqo}).
A {\it super-word} is a word in super-letters. 
\label{sup1}
\end{definition}

By Shirshov theorem every standard word $u$
defines the only super-letter, in what follows we will denote it by $[u].$
The order on the super-letters is defined in the natural way:
$[u]>[v]\iff u>v.$
We should stress that this order is not complete: 
If $x>y,$ there are infinite chains of super-letters
$$
[xy]>[xy^2]>\ldots >[xy^n]>\ldots \ ; [xy]<[x^2y]<\ldots <[x^ny]<\ldots \, .
$$
Nevertheless the Hall order on super-words is still a complete monoidal order. 
\begin{definition} \rm
A super-letter $[u]$
is called {\it hard in }$H$
provided that its value in $H$
is not a linear combination
of values of super-words in smaller than $[u]$ super-letters.
\label{tv1}
\end{definition}
\begin{definition} \rm
We say that a {\it height} of a hard in $H$ super-letter $[u]$
equals $h=h([u])$ if  $h$
is the smallest number such that: first, $p_{uu}$
is a primitive $t$-th root of unity and either $h=t$
or $h=tl^r,$ where $l=$char({\bf k}); and  then the value in $H$
of $[u]^h$
is a linear combination of super-words in  less than $[u]$ super-letters.
If there exists no such number  then the height equals infinity.
\label{h1}
\end{definition}
\begin{theorem}
Values of all hard in $H$ super-letters form a set of PBW-generators of $H$
over {\bf k}$[G].$
\label{BW}
\end{theorem}
\begin{proof} If $X$ is a finite set, this statement has
been proved in \cite[Theorem 2]{Kh3}, see a footnote 
on page 267. In the general case the algebra
$G\langle X\rangle $ has a covering by 
$G\langle X_{\alpha }\rangle ,$ where $X_{\alpha }$
runs through all finite subsets of $X.$ We shall consider each of
$X_{\alpha }$ as an ordered subset of $X,$
and $G\langle X_{\alpha }\rangle $ as a $\Gamma $-graded subalgebra
of $G\langle X\rangle .$ Respectively $H$ has a covering
by the finitely generated $\Gamma $-graded character Hopf subalgebras
$H_{\alpha }\stackrel{df}{=}\xi (G\langle X_{\alpha }\rangle ).$

By Definition \ref{tv1} every hard in $H$
super-letter $[u]$ is hard in $H_{\alpha }$
as soon as $u\in G\langle X_{\alpha }\rangle .$
If a super-letter $[v]\in G\langle X_{\alpha }\rangle$ is not hard in $H,$
then by definition $\xi ([v])=\xi (A),$ 
where $A$ is a linear
combination of super-words in less than $[v]$ super-letters.
Since  $H$ is $\Gamma $-graded,  all super-words
of the sum $A$ have the same constitution as the word $v$ does.
In particular all super-words of $A$
belong to $G\langle X_{\alpha }\rangle .$ Hence $[v]$ is
not hard in $H_{\alpha }.$

In the same way the height function is independent of $\alpha .$
\end{proof}
 
The following statements provide important properties of 
PBW-basis defined by the hard super-letters.
\begin{lemma}$(${\rm \cite[Lemma 8]{Kh3}}$).$
The coproduct of a super-letter $[w]$ has a representation
\begin{equation}
\Delta ([w])=[w]\otimes 1+g_w\otimes [w]+
\sum _i\alpha _ig(W_i^{\prime \prime })W_i^{\prime }
\otimes W_i^{\prime \prime },
\label{[w]}
\end{equation}
where $W_i^{\prime }$ are non-empty words in less than $[w]$
super-letters.
\label{ume}
\end{lemma}
\begin{lemma} The Hall order of the super-words is compatible with
the PBW- decomposition related to the hard super-letters
$($see Definition $\ref{pbgd}).$
\label{redu}
\end{lemma}
\begin{proof} 
Let $W$ be a super-word.
There exists the following natural diminishing process of the decomposition. First,
in {\bf k}$\langle X\rangle $ according to
\cite[Lemma 7]{Kh3} we decompose the super-word $W$ 
in a linear combination of smaller monotonous super-words, then, we replace  
each non hard super-letter with the decomposition of its value in $H$ that exists 
by Definition \ref{tv1}, and again decompose the appeared super-words 
in linear combinations of smaller monotonous super-words, and so on, until we get a linear combination of monotonous super-words in hard
super-letters. If they are not restricted, we may apply Definition \ref{h1} and repeat the process
until we get only monotonous restricted words in hard super-letters.\end{proof}
\begin{lemma}
If $[w]$ is a super-letter then 
\begin{equation}
\Delta (\xi ([w]^{m}))=\sum _{j=0}^{m}\Bigl[^m_j\Bigr]_q g_w^{m-j}\xi ([w])^j\otimes 
\xi ([w])^{m-j}+
\sum _i\alpha _ig(V_i)\xi (U_i)\otimes \xi (V_i), 
\label{[w]m}
\end{equation}
where $\Bigl[^m_j\Bigr]_q$ are the $q$-binomial coefficients 
with $q=p(w,w),$ and $U_i$ are basis super-words
such that the first from the left super-letter is less than $[w].$ 
\label{slep}
\end{lemma}
\begin{proof} The $m$-th power of the right hand side 
of (\ref{[w]}) after developing of the product
takes the form (\ref{[w]m}), where each of $U_i$
is a product of $m$ super-words some of whom equal to 
$[w]$ (but not all of them!) and others equal to $W_i^{\prime }$'s.
Hence $U_i$ as a super-word is less than $[w]^{m}$
with respect to the lexicographical ordering of words in super-letters.
Let us decompose $\xi (U_i)$ in the PBW-basis defined by the hard super-letters.
According to Lemma 
\ref{redu} we turn to the formula (\ref{[w]m})
where still $U_i<[w]^m.$ 
This implies the required property since the first super-letter
in  a basis super-word is always the minimal one.
\end{proof}

\section{Right coideal subalgebras}
\begin{theorem} Let $H$ be a character Hopf algebra. Every  
right coideal subalgebra {\bf U}$\supseteq {\bf k}[G]$ has a PBW-basis
that may be extended up to a PBW-basis of $H.$
\label{max}
\end{theorem}
By Lemma \ref{gruf} it suffices to show that gr{\bf U} as a subalgebra of 
gr$H$ has a required basis. Since gr{\bf U} is a right coideal subalgebra of gr$H,$
we still may suppose that $H$ is a $\Gamma $-graded character Hopf algebra
and {\bf U} is a $\Gamma $-graded coideal subalgebra.
\begin{lemma} 
Let $w\in {\bf U}$.  If $\pi , \mu :H\rightarrow H$ are linear maps such that 
$(\pi \otimes \mu )\Delta (w)=a\otimes b$
with $b\neq 0,$ then $a=\pi (c)$ for some  $c\in {\bf U}.$
\label{tool}
\end{lemma}
\begin{proof} We have
$(\pi \otimes \mu )\Delta ({\bf U})\subseteq \pi ({\bf U})\otimes H.$
Hence $a\otimes b\in \pi ({\bf U})\otimes H,$ and $a\in \pi ({\bf U}).$
\end{proof}

We shall use this evident statement as a basic tool. Note that once we have a PBW-basis
of $H,$ to define a linear map it suffices to fix its  values on the restricted
monotonous words in an arbitrary way.

Let $[u]$ be a hard super-letter.
Suppose that in {\bf U} there exists an element $c$ with the leading super-word $[u]^m,$
$m\geq 1.$ Since $G\subseteq ${\bf U},
we may suppose that the super-word $[u]^m$ appears one time with the trivial coefficient:  
\begin{equation}
c=\xi ([u]^m+\sum _j\alpha _jg_j[u]^m+\sum _i\alpha _ig_iW_iR_i), 
\label{khu}
\end{equation}
where  $W_{i}$ are  nonempty basis super-words  in less than $[u]$ super-letters, 
while $R_i$ are basis super-words in greater than or equal to $[u]$ super-letters,
$\alpha _{i}, \alpha _{j}\in {\bf k}, g_{i}, g_{j}\in G,$ $g_j\neq 1.$

Denote by $\iota $ the natural projection $H\rightarrow \, ${\bf k}, $\iota (g\xi (W))=0,$
unless $g=1,$ $W=\emptyset .$ 
Since $\Delta (gW)=(g\otimes g)\Delta (W),$ we have
$$
({\rm id} \otimes \iota )(\Delta (c))=\xi ([u]^m+\sum _{g_i=1}\alpha _{i}W_{i}R_i)
\otimes 1.
$$
Thus, by Lemma \ref{tool}
\begin{equation}
c^{\prime }=\xi ([u]^m+\sum _{g_i=1}\alpha _i W_iR_i)\in {\bf U}.
\label{1c}
\end{equation}
\underline{In what follows we fix the notation} $c_u$ for one of the elements from
{\bf U} of the form $($\ref{1c}$)$ that has the minimal possible $m.$ 

\begin{lemma} In the representation $(\ref{1c})$ of the chosen element $c_u$
either $m=1,$ or $p_{uu}$ is a primitive $t$-th root of unity and $m=t$ or 
$($in the case of positive characteristic$)$ $m=t({\rm char}\, {\bf k})^s.$
In particular $c_u$ is a thin element, see $(\ref{thin}).$ 
\label{co1}
\end{lemma}
\begin{proof}
If $m=1$ there is nothing to prove. Let $m>1.$
For each $k, 1\leq k<m$ we consider the following  linear map
set up on the PBW-basis of super-words.
\begin{equation}
\pi _k(g\xi (W))=\left\{ \begin{matrix}
0, &\ \ \hbox{if } \ \ W\prec [u]^k,\hfill \cr
g\xi (W), & \ \ \hbox{otherwise;}\hfill 
 \end{matrix} \right. 
\label{2c}
\end{equation}
By means of formula (\ref{[w]m}) we have 
\begin{equation}
(\pi _k\otimes \pi _{m-k})\Delta (\xi ([u])^m)=\Bigl[ ^m_k\Bigr]_q g_u^{m-k}
\, \xi ([u])^k\otimes \xi ([u])^{m-k}.
\label{4c}
\end{equation}
By Lemma \ref{ume} and Lemma \ref{redu} the coproduct $\Delta (\xi (W_i))$
is the sum of tensors  $g\xi(W_i^{\prime })\otimes \xi(W_i^{\prime \prime}),$
where all $W_i^{\prime }$ are the basis super-words lexicographically smaller than
$[u]^m$ with the only exception that equals $g(W_i)\otimes \xi (W_i).$
Hence $(\pi_k\otimes \pi _{m-k})\Delta (\xi (W_iR_i))$ is a sum of tensors 
of the form $g\pi _k(\xi(R^{\prime}_i))\otimes \pi _{m-k}(\xi (W_iR_i^{\prime \prime})).$
Again by Lemma \ref{redu} we have $\pi _{m-k}(\xi (W_iR_i^{\prime \prime}))=0.$ Hence
\begin{equation}
(\pi _k\otimes \pi _{m-k})\Delta (\xi (W_iR_i))=0.
\label{5c}
\end{equation}
Thus we may write
\begin{equation}
(\pi _k\otimes \pi _{m-k})\Delta (c_u)=
\Bigl[ ^m_k\Bigr]_q g_u^{m-k}
\, \xi ([u])^k\otimes \xi ([u])^{m-k}.
\label{7c}
\end{equation}
If $\Bigl[ ^m_k\Bigr]_q \neq 0,$ then by Lemma \ref{tool}
we find $c\in {\bf U}$ such that $\pi _k(c)=g_u^{m-k}\xi ([u]^k).$
By definition of $\pi _k$ this equality means that the 
PBW-decomposition of $g_u^{k-m}c$
has the form (\ref{1c}) with $k$ in place of $m.$ 
Since by the choice of $c_u$ the number $m$ is minimal,
we get $\Bigl[ ^m_k\Bigr]_q=0,$ $1\leq k<m.$
This system of equations implies $q^m=1$ and either $m$ equals
the multiplicative order $t$ of $q$ or $m=t({\rm char }\, {\bf k})^s.$
\end{proof}

By Lemma \ref{pbw} the set $T$ of all above defined elements $c_u$
has an extention up to a set $P_T$  of PBW-generators of $H$ over ${\bf k}[G].$
Now Theorem \ref{max} follows from the proposition below.
\begin{proposition}  An element $c\in H$ belongs to
{\bf U} if and only if all PBW-generators in the PBW decomposition of $c$
with respect to $P_T$ belong to $T.$ In particular $T$ is a set of PBW-generators
of {\bf U} over ${\bf k}[G].$
\label{main}
\end{proposition}
To prove this statement we will need some additional properties of the PBW-basis
defined by $P_T.$ We extend the order ``$<$" already defined on $P_T$ onto 
the set of all words in $P_T$ as the lexicographical order. 
The order ``$\prec $" 
is the Hall order induced by the degree function $D;$ that is, $W\prec U$
if and only if either $D(W)<D(U),$ or $D(W)=D(U)$ $\& $ $W<U.$
We have to stress that the order on the set of letters $P_T$ differs from the
order on the set of one-letter words $P_T.$ 

We start with connections between these two PBW-decompositions.
\begin{lemma} The leading term of a monotonous restricted 
word in $P_T,$
\begin{equation}
\theta _1^{n_1}[u_1]^{r_1}\theta _2^{n_2}[u_2]^{r_2}\cdots \theta _k^{n_k}[u_k]^{r_k}, 
\label{slot}
\end{equation}
under the PBW-decomposition related to the hard super-letters equals
\begin{equation}
[u_1]^{n_1m_1+r_1}[u_2]^{n_1m_2+r_2}\cdots [u_k]^{n_km_k+r_k},
\label{slo1t}
\end{equation}
where $\theta _i=[u_i]^{m_i}+\cdots .$ Conversely if 
\begin{equation}
W=[u_1]^{s_1}[u_2]^{s_2}\cdots [u_k]^{s_k}
\label{slot2}
\end{equation}
is a monotonous restricted super-word in hard super-letters,
 then its leading term in the decomposition
with respect to $P_T$ equals $(\ref{slot}),$ where $n_i=\lfloor s_i/m_i\rfloor ,$
$r_i=s_i-n_im_i.$
\label{pbz}
\end{lemma}
\begin{proof} The first part of the lemma has been proved in Lemma \ref{pbw}.
The second part follows by induction on super-words ordered by the Hall order.
Indeed the difference $E$ between $W$ and (\ref{slot}) is a linear combination
of super-words that are less than $W.$ By Lemma \ref{redu} all basis super-words
in the PBW-decomposition of $E$ are less than $W.$ It remains to note
that if $W^{\prime }\prec W$ is another basis super-word, then the word
(\ref{slot}) related to $W^{\prime }$ is less than that related to $W.$
\end{proof}
\begin{lemma}
The Hall order on the words in $P_T$ is compatible with the 
PBW-decomposition related to $P_T$ $($see Definition $\ref{pbgd}).$
\label{redz}
\end{lemma}
\begin{proof} Let $W$ be a word in $P_T.$
The word $W$ has the form (\ref{slot}),
where  $\theta _i=c_{u_i}=\xi ([u_i]^{m_i}+\cdots ),$ while  $u_i$'s are not necessary
increase from the left to the right.  
If we replace $\theta _i:=[u_i]^{m_i}+\cdots $
in (\ref{slot}) and then develop the multiplication in $G\langle X\rangle ,$
we get a linear combination $\Sigma $ over {\bf k}$[G]$ of super-words  with the leading term 
\begin{equation}
[u_1]^{n_1m_1+r_1}[u_2]^{n_1m_2+r_2}\cdots [u_k]^{n_km_k+r_k}.
\label{slot3}
\end{equation}
By Lemma \ref{redu} each super-word $W_1$ in the PBW-decomposition
 of $\xi (\Sigma )$ is less than or equal to (\ref{slot3}). That is, 
if  $W_1$ is different from (\ref{slot3}), we may write
\begin{equation}
W_1=[u_1]^{n_1m_1+r_1}[u_2]^{n_2m_2+r_2}\cdots [u_s]^{n_sm_s+r_s}[u_{s+1}]^{t}[v]\cdots ,
\label{slo2t}
\end{equation}
here $[v]$ is the first from the left super-letter where $W_1$ differ from (\ref{slot3}),
and $0\leq t<n_{s+1}m_{s+1}+r_{s+1}.$ Since $W_1$ is monotonous, we have
$v>u_{s+1},$ provided that $t\neq 0.$ However in this case $W_1$ is greater than (\ref{slot3}).
Hence $t=0,$ and $u_{s+1}>v\geq u_s\geq u_{s-1}\geq \ldots \geq u_1.$ 

If all inequalities among $u_i$'s in the above chain are strict,
then according to Lemma \ref{pbz}, all words
of the PBW-decomposition with respect to $P_T$ of $\xi (W_1)$ are less than or equal to
\begin{equation}
\theta _1^{n_1}[u_1]^{r_1}\theta _2^{n_2}[u_2]^{r_2}\cdots \theta _s^{n_s}[u_s]^{r_s}\theta _v, 
\label{slo3}
\end{equation}
where $\theta _v=c_v,$ or $\theta _v=[v].$ Hence they are less than $W.$ 

If not all inequalities are strict, say $u_1<u_2<\ldots <u_p=\cdots =u_q<u_{q+1},$ $q\leq s,$ then 
again by Lemma \ref{pbz} the $P_T$-leading word of $\xi (W_1)$ starts with
\begin{equation}
\theta _1^{n_1}[u_1]^{r_1}\cdots \theta _{p-1}^{n_{p-1}}[u_{p-1}]^{r_{p-1}}
\theta _p^{(n_p+\cdots +n_q+n_0)m_p}[u_p]^{r_0}, 
\label{slo4}
\end{equation}
where $r_p+\cdots +r_q=n_0m_p+r_0,$ $0\leq r_0<m_p.$
It is still less than $W$ since by definition $\theta _p <[u_p],$ provided 
that $ m_p>1.$ \end{proof}

\begin{lemma}
The coproduct of each ${\theta }\in P_T$ has a representation
\begin{equation}
\Delta (\theta )=\theta \otimes 1+g_{\theta }\otimes \theta +
\sum _i\alpha _ig_iW_i^{\prime }
\otimes W_i^{\prime \prime },
\label{cobuk}
\end{equation}
where $W_i^{\prime },$ $W_i^{\prime \prime }$ are restricted monotonous words $($products$)$
 in $P_T,$
and for every $i$ either $W_i^{\prime }$ or $W_i^{\prime \prime }$ starts with a letter that is less than
$\theta .$ 
\label{ume1}
\end{lemma}
\begin{proof} 
Let $\theta =\xi ([u]^m+\sum _iW_iR_i)$ be the decomposition
of $\theta \in P_T$ with respect to the hard super-letters. By Lemma \ref{slep} and Lemma \ref{co1},
we have
$$
\Delta (\xi ([u]^m))-\xi ([u]^m)\otimes 1-g_u^m\otimes \xi ([u]^m)=
\sum _i\alpha _ig_i\xi(U_i)\otimes \xi (V_i), 
$$
where $U_i$ are basis super-words
 starting with smaller than $[u]$ super-letters. By the second part 
of Lemma \ref{pbz} each summand of the decomposition of $\xi (U_i)$
in the PBW-basis defined by $P_T$ starts by $c _w$ or $[w]$ with  $w<u.$

Let $W_i=[w]U,$ where $[w]<[u].$ By Lemma \ref{ume} we have
$$
\Delta (W_iR_i)=\sum_j (V_{ij}^{\prime }\otimes V_{ij}^{\prime \prime })\Delta (UR_i)
+(g_w\otimes [w])\Delta (UR_i),
$$ 
where all $V_{ij}^{\prime }$ are non empty super-words in less than $[u]$ super-letters.
All left hand sides of the tensors  
$(V_{ij}^{\prime }\otimes V_{ij}^{\prime \prime })\Delta (UR_i)$ start with smaller than $[u]$
super-letters.
Hence by Lemma \ref{redu} it remains to use the converse part of Lemma \ref{pbz}.

In perfect analogy, if
$\Delta (UR_i)=\sum _jU_{ij}^{\prime }\otimes U_{ij}^{\prime \prime },$
then $(g_w\otimes [w])\Delta (UR_i)=g_wU_{ij}^{\prime }\otimes [w]U_{ij}^{\prime \prime }R_i.$
Therefore the right hand side
of each tensor  that appears 
in the $P_T$-decomposition of $(g_w\otimes \xi ( [w]))\Delta (\xi (UR_i))$ 
starts with a letter that is less than $\theta .$
\end{proof}
\begin{lemma}
Let ${\theta }\in P_T.$ The coproduct of $\theta ^n$ has a decomposition
\begin{equation}
\Delta (\theta ^n)=\sum _{j=0}^n\Bigl[^n_j\Bigr]_q g_{\theta }^{n-j}\theta ^j
\otimes \theta ^{n-j} +\sum _i\alpha _ig_iW_i^{\prime } \otimes W_i^{\prime \prime },
\label{ncobuk}
\end{equation}
where 
$\Bigl[^n_j\Bigr]_q$ are $q$-binomial coefficients with $q=p(\theta ,\theta ),$
$W_i^{\prime },$ $W_i^{\prime \prime }$ are restricted monotonous words in $P_T,$
and for every $i$ either $W_i^{\prime }$ or $W_i^{\prime \prime }$ starts with a letter that is less than
$\theta .$ 
\label{ume2}
\end{lemma}
\begin{proof} 
If we develop multiplication in the $n$-th power of the right hand side of (\ref{cobuk}), 
then we get the first sum of (\ref{ncobuk}) and a  {\bf k}$[G]$-linear combination of tensors 
$\theta ^iW^{\prime }\otimes \theta ^jW^{\prime \prime},$
where $i$ and $j$ may be zero, but either $W^{\prime }$ or $W^{\prime \prime}$ starts with
a less than $\theta $ letter. Let it be $W^{\prime }.$ 
By Lemma \ref{redz} the PBW-decomposition
of $\theta ^iW^{\prime }$ has only words that are less than $\theta ^iW^{\prime }$ 
in lexicographical order. Since these words are monotonous,
no one of them may start with a letter that 
is grater than or equal to $\theta .$ The lemma is proved.\end{proof}

 {\it Proof }of Proposition \ref{main}. Let in contrary
$c\in {\bf U}$ be an element of the minimal degree whose decomposition
in $P_T$-basis has super-letters $[u]\in P_T\setminus T.$
Since {\bf U} is a subalgebra, we may suppose that each term of the 
decomposition has such a super-letter.

Let $U=\theta _1^{n_1}\theta _2^{n_2}\cdots \theta _k^{n_k},$ 
$ \theta_1<\theta_2<\cdots <\theta_k$ be the leading word
of $c$ in the $P_T$-basis.
Since $G\subseteq ${\bf U},
we may suppose that $U$ appears one time with the trivial coefficient:  
\begin{equation}
c=U+\sum _j\alpha _jg_jU+\sum _i\beta _ig_iU_i, 
\label{mkhu}
\end{equation}
where
$D(U_{i})=D(U),$
$\alpha _{j},$ $\beta_{i}\in {\bf k},$  $g_{i}, g_j\in G,$ $g_j\neq 1.$

Denote by $\iota $ the natural projection $H\rightarrow \, ${\bf k}, $\iota (gW)=0,$
unless $g=1,$ $W=\emptyset .$ 
Here $W$ is an arbitrary restricted monotonous word in $P_T.$
Since $\Delta (gW)=(g\otimes g)\Delta (W),$ we have
$$
({\rm id} \otimes \iota )(\Delta (c))=(U+\sum _{g_i=1}\beta _{i}U_i)\otimes 1,
$$
Thus, by Lemma \ref{tool}
\begin{equation}
c^{\prime }=U+\sum _{g_i=1}\beta _i U_i\in {\bf U} 
\label{m1c}
\end{equation}
where $U_i<U.$ To get a contradiction we consider two cases.

{\it Case} 1. Suppose that $\theta _k\notin T;$ that is, $\theta _k=[u]$ and in {\bf U} there does not
exist an element with the leading super-word $[u]^m,$ $1\leq m\leq n_k.$
Let us define the following two linear maps
set up on the PBW-basis related to $P_T.$
\begin{equation}
\pi (gW)=\left\{ \begin{matrix}0, &\  \hbox{if }\ \ W\prec \theta _k^{n_k}; \hfill \cr
gW, & \ \hbox{otherwise.}\hfill 
\end{matrix} \right. 
\label{m2c}
\end{equation}
\begin{equation}
\nu (gW)=\left\{ \begin{matrix}0, &\ \hbox{if }\ \
W\prec \theta _1^{n_1}\cdots \theta _{k-1}^{n_{k-1}}; \hfill \cr
gW, & \ \hbox{otherwise.}\hfill 
\end{matrix}\right.
\label{m3c}
\end{equation}
Let us show first that 
\begin{equation}
(\pi \otimes \nu )(\Delta (U))=g\theta _k^{n_k}\otimes \theta _1^{n_1}\cdots \theta _{k-1}^{n_{k-1}},
\ \ \ g\in G.
\label{mac}
\end{equation}
Since the word $U$ is monotonous, by Lemma \ref{ume2} we have
\begin{equation}
\Delta (\theta _i^{n_i})=
\sum_{j_i} \alpha _{j_i}g_{j_i}W^{\prime }_{j_i}\otimes W^{\prime \prime }_{j_i},\ \ i<k,
\label{mac1}
\end{equation}
where for each $j_i,$ $i\leq k$ either  $W^{\prime }_{j_i}$ 
starts by a letter from $P_T$ that is less than $\theta _k,$
or  $W^{\prime \prime }_{j_i}$
starts by a letter from $P_T$ that is less than $\theta _i ,$
with the only exception  $W^{\prime }_{j_i}=\emptyset ,$ $W^{\prime \prime }_{j_i}=\theta _i^{n_i}.$
In the same way
\begin{equation}
\Delta (\theta _k^{n_k})=
\sum_{j_k}  \alpha _{j_k}g_{j_k}W^{\prime }_{j_k}\otimes W^{\prime \prime }_{j_k}, 
\label{mac2}
\end{equation}
where either  $W^{\prime }_{j_k}$ 
or  $W^{\prime \prime }_{j_i}$
starts by a letter from $P_T$ that is less than $\theta _k,$
with  exceptions of the form   $W^{\prime }_{j_k}$ $=\theta _k^j ,$ $W^{\prime \prime }_{j_k}$ 
$=\theta _k^{n_k-j}.$

$\Delta (U)$ is a right 
linear combination over ${\bf k}[G]$ of tensors
\begin{equation}
E=W^{\prime }_{j_1}W^{\prime }_{j_2}\cdots W^{\prime }_{j_k}\otimes 
W^{\prime \prime }_{j_1}W^{\prime \prime }_{j_2}\cdots W^{\prime \prime }_{j_k}.
\label{mac3}
\end{equation}
Let $W^{\prime }_{j_s}$ be the first from the left nonempty factor in $E.$ 

If, first, 
$s\neq k,$ then the tensor get the form
$$
E=W^{\prime }_{j_s}\cdots W^{\prime }_{j_k}\otimes 
\theta _1^{n_1}\cdots \theta _{s-1}^{n_{s-1}}
W^{\prime \prime }_{j_s}\cdots W^{\prime \prime }_{j_k}.
$$

If $W^{\prime }_{j_s}$ starts by a less than $\theta _k $ letter, 
then the left hand side of $E$ is less than $\theta _k^{n_k},$
hence by Lemma \ref{redz} we get $(\pi \otimes \nu)(E)=0.$

If $W^{\prime \prime }_{j_s}$ starts by a less than $\theta _s$ letter,
then the right hand side is less than $\theta _1^{n_1}\cdots \theta _{k-1}^{n_{k-1}},$
and again $(\pi \otimes \nu)(E)=0.$

If, next, $s=k,$ then the tensor has the form 
$$
E=W^{\prime }_{j_k}\otimes  \theta _1^{n_1}\cdots \theta _{k-1}^{n_{k-1}}
W^{\prime \prime }_{j_k}.
$$
Since $D(W^{\prime }_{j_k})=n_kD(\theta _k)-D(W^{\prime \prime }_{j_k})<n_kD(\theta _k),$
unless $W^{\prime \prime }_{j_k}=\emptyset ,$
we get $\pi (W^{\prime }_{j_k})=0$ with the only exception
$W^{\prime }_{j_k}=\theta _k^{n_k}.$ 
Thus $\pi \otimes \nu$ kills all tensors in the PBW-decomposition of $\Delta (U)$
except one of them. This proves (\ref{mac}).

Let us show, further, that
\begin{equation}
(\pi \otimes \nu)(U_i)=0.
\label{macr}
\end{equation}
Since $U_i<U,$ we have 
$U_i=\theta _1^{n_1}\cdots \theta _{s}^{n_s}\theta _{s+1}^r\eta \cdots ,$
where  $0\leq r<n_{s+1},$ $\eta <\theta _{s+1}.$ If $r\neq 0,$
then  $\theta _{s+1}<\eta $ since the word $U_i$ is monotonous. Thus $r_s=0;$
that is, 
$$
U_i=\theta _1^{n_1}\cdots \theta _{s}^{n_s}\eta \cdots , \ \ \ \eta <\theta _{s+1}.
$$
By means of Lemma \ref{ume1} we have
\begin{equation}
\Delta (\eta )=
\sum_j \alpha _jg_jV^{\prime }_j\otimes V^{\prime \prime }_j,
\label{macm}
\end{equation}
where for each $j,$ either  $V^{\prime }_j$ or $V^{\prime \prime }_j$
starts by a letter from $P_T$ which is less than $\theta _{s+1}.$
Thus $\Delta (U_i)$ is a right 
linear combination over ${\bf k}[G]$ of tensors
\begin{equation}
E=W^{\prime }_{j_1}\cdots W^{\prime }_{j_s}V_j^{\prime }\cdots \otimes 
W^{\prime \prime }_{j_1}\cdots W^{\prime \prime }_{j_s}V_j^{\prime \prime }\cdots .
\label{mac3m}
\end{equation}

If  $W^{\prime }_{j_1}=\cdots = W^{\prime }_{j_s}=V_j^{\prime }=\emptyset ,$
then the right hand side takes up the form
 $\theta _1^{n_1}\cdots \theta _{s}^{n_s}\eta \cdots $ which is less than
$\theta _1^{n_1}\cdots \theta _{k-1}^{n_{k-1}},$ hence $(\pi \otimes \nu)(E)=0.$

Let $W^{\prime }_{j_r},$ $r\leq s$ be the first from the left nonempty word in $E.$ 
The tensor takes up the form
$$
E=W^{\prime }_{j_r}\cdots W^{\prime }_{j_s}V_j^{\prime }\cdots \otimes 
\theta _1^{n_1}\cdots \theta _{r-1}^{n_{r-1}}
W^{\prime \prime }_{j_r}\cdots W^{\prime \prime }_{j_s}V_j^{\prime \prime }\cdots.
$$
Since either $W^{\prime }_{j_r}$ starts with a less than $\theta _k $ letter, 
or $W^{\prime \prime }_{j_r}$ starts with a less than $\theta _r$
letter, again $(\pi \otimes \nu)(E)=0.$

If the first from the left nonempty factor is $V_j^{\prime },$ then
in the same way either $V_j^{\prime }$ starts with a letter $\leq \eta <\theta _{s+1}\leq \theta _k,$
or $V_j^{\prime \prime }$ does. In both cases $(\pi \otimes \nu)(E)=0.$

Thus
$(\pi \otimes \nu)(\Delta (U_i))=0.$ Now taking into account (\ref{m1c}),
we get 
\begin{equation}
(\pi \otimes \nu )(\Delta (c^{\prime }))=
g\theta _k^{n_k}\otimes \theta _1^{n_1}\cdots \theta _{k-1}^{n_{k-1}}, \ g\in G.
\label{ma2c}
\end{equation}
Lemma \ref{tool} implies that in {\bf U} there exists an element $c^{\prime \prime }$
such that $\pi (c^{\prime \prime })=g\theta _k^{n_k}.$ Definition (\ref{m2c}) of $\pi $
shows that $\theta _k^{n_k}=[u]^{n_k}$ is the leading word of $c^{\prime \prime }$
in the $P_T$-decomposition. By Lemma \ref{pbz} the leading super-word 
of $c^{\prime \prime }$ equals $[u]^{n_k}.$ According to the construction
of $T,$ there exists an element $c_u\in T$ of the form
$c_u=[u]^m+\cdots ,$ $m\leq n_k.$ This contradicts the conditions of the first case
$([u]=\theta _k\notin T, h_T([u])>n_k).$

{\it Case} 2. Suppose  that $\theta _k\in T.$ In this case $\theta _k\in {\bf U},$
therefore one of the letters $\theta _1, \ldots , \theta _{k-1}$ does not belong to $T.$
By the inductive supposition no one element with the leading 
word $\theta _1^{n_1}\cdots \theta _{k-1}^{n_{k-1}}$ belongs to {\bf U}.
At the same time, in perfect analogy with the first case (up to left-right
symmetry in the consideration of tensors), we have 
$$
(\nu \otimes \pi )(\Delta (c^{\prime }))=g\theta _1^{n_1}\cdots \theta _{k-1}^{n_{k-1}}\otimes \theta _k^{n_k}.
$$
Hence by Lemma \ref{tool} there exists $c^{\prime \prime }\in {\bf U}$ such that 
$\nu (c^{\prime \prime })=g\theta _1^{n_1}\cdots \theta _{k-1}^{n_{k-1}}.$
Definition (\ref{m3c}) of $\nu $ implies that the leading word of $c^{\prime \prime }$
indeed equals $\theta _1^{n_1}\cdots \theta _{k-1}^{n_{k-1}}.$
A contradiction. Proposition \ref{main} and Theorem \ref{max} are completely proved.

\bibliographystyle{amsplain}

\end{document}